\documentclass[11pt,a4paper,reqno]{amsart}
\usepackage{amsmath,amsfonts,amssymb,amsthm,mathrsfs,mathtools}
\usepackage{microtype,enumitem,comment}
\usepackage[dvipsnames]{xcolor}
\usepackage[pagebackref=true,colorlinks=true]{hyperref}
\hypersetup{
  linkcolor=blue,
  citecolor=blue,
  urlcolor=blue,
  filecolor=blue,
  pdfborder={0 0 0},
  pdftitle={Strong endpoint estimates for bilinear spherical averaging operators on closed hypersurfaces and curved patches},
  pdfauthor={Xi Cen}
}

\renewcommand*{\backref}[1]{}

\renewcommand*{\backrefalt}[4]{%
  \ifnum#1>0\relax
    \space{\footnotesize\color{blue}#2}%
  \fi
}

\numberwithin{equation}{section}
\theoremstyle{plain}
\newtheorem{theorem}{Theorem}[section]
\newtheorem{problem}{Problem}
\newtheorem{lemma}[theorem]{Lemma}

\newtheorem{corollary}[theorem]{Corollary}
\theoremstyle{definition}
\newtheorem{remark}[theorem]{Remark}
\newcommand{\R}{\mathbb R}
\newcommand{\Sph}{\mathbb S}
\newcommand{\A}{\mathcal A}
\newcommand{\B}{\mathcal B}
\newcommand{\pn}{\pi_{-}}
\newcommand{\dd}{\,\mathrm d}
\newcommand{\one}{\mathbf 1}
\newcommand{\spt}{\operatorname{supp}}
\newcommand{\norm}[1]{\lVert#1\rVert}
\newcommand{\arxiv}[1]{\href{https://arxiv.org/abs/#1}{arXiv:\nolinkurl{#1}}}
\DeclareRobustCommand{\Thmref}[1]{\mbox{\hyperref[#1]{Theorem~\ref*{#1}}}}
\DeclareRobustCommand{\Probref}[1]{\mbox{\hyperref[#1]{Problem~\ref*{#1}}}}
\DeclareRobustCommand{\Lemref}[1]{\mbox{\hyperref[#1]{Lemma~\ref*{#1}}}}
\DeclareRobustCommand{\Propref}[1]{\mbox{\hyperref[#1]{Proposition~\ref*{#1}}}}
\DeclareRobustCommand{\Corref}[1]{\mbox{\hyperref[#1]{Corollary~\ref*{#1}}}}
\DeclareRobustCommand{\Remref}[1]{\mbox{\hyperref[#1]{Remark~\ref*{#1}}}}
\DeclareRobustCommand{\Secref}[1]{\mbox{\hyperref[#1]{Section~\ref*{#1}}}}
\title[Strong endpoints for bilinear spherical averages]{\bf Strong endpoint estimates for bilinear spherical averaging operators on closed hypersurfaces and curved patches}
\author[X. Cen]{Xi Cen}
\address{Xi Cen\\
School of Science\\
China University of Mining and Technology-Beijing\\
Beijing 100083, People's Republic of China}
\email{xicenmath@gmail.com}
\date{}
\begin{document}
\begin{abstract}
Tainara Borges~\cite[Section~7.1, Question~7.2]{Borges2026}
asked whether the strong $L^1\times L^1\to L^{1/2}$ estimate for
bilinear spherical averages extends to all compact smooth
hypersurfaces without boundary and with nonvanishing Gaussian curvature.
We give a positive answer when $d\ge2$. The geometric argument
concerns the difference projection $(y,z)\mapsto y-z$: at every
critical point, two definite directions in its fiber already ensure
that the projected surface measure has bounded density. This
condition holds for the closed hypersurfaces in the original question.
The ensuing strong endpoint uses the localization principle of
Iosevich, Palsson and Sovine~\cite[Propositions~4.1--4.2]{IPS2022}.
Our proof also establishes a partial positive result for the
curvature-rank problem in Borges~\cite[Question~7.1]{Borges2026}.
We exhibit a smooth full-curvature patch with boundary for which
the strong estimate fails. The endpoint yields the entire fixed-scale
H\"older family.
\end{abstract}
\maketitle
\thispagestyle{plain}

\section{Introduction}\label{sec:intro}
\subsection{Background and motivations}\label{subsec:background}
For $d\ge1$ and $t>0$, the bilinear spherical averaging operator is
\begin{equation}\label{eq:spherical}
 \A_t(f,g)(x)\coloneqq\int_{\Sph^{2d-1}}
 f(x-ty)g(x-tz)\dd\sigma(y,z),\qquad x\in\R^d,
\end{equation}
where $\sigma$ is normalized surface measure. The two translations are tied to a single point on the sphere in $\R^{2d}$. This kind of multilinear convolution was studied by Oberlin~\cite{Oberlin1988}. Jeong and Lee~\cite{JeongLee2020} considered the associated maximal averages and bilinear Bochner--Riesz operators.

For the fixed-scale operator, the first question is what happens when both
functions are merely integrable. The answer on the sphere is known. When
$d\ge2$, Iosevich, Palsson and Sovine~\cite[Section~5.2]{IPS2022} proved
\begin{equation}\label{eq:knownsphere}
 \norm{\A_1(f,g)}_{L^{1/2}(\R^d)}
 \lesssim_d\norm{f}_{L^1(\R^d)}\norm{g}_{L^1(\R^d)}.
\end{equation}
The relevant geometry is seen after integration in $x$. For a positive compactly supported measure $\mu$ on $\R^{2d}$, set
$T_\mu(f,g)(x)=\int f(x-y)g(x-z)\dd\mu(y,z)$ and
$\nu=(y-z)_*\mu$. For nonnegative $f,g$,
\begin{equation}\label{eq:motivationidentity}
 \int_{\R^d}T_\mu(f,g)(x)\dd x
 =\int_{\R^d}\!\int_{\R^d} f(u)g(u+v)\dd u\dd\nu(v).
\end{equation}
In particular, a bounded density for $\nu$ gives an $L^1$ estimate.
The extra passage to strong $L^{1/2}$ is supplied by the localization
principle of Iosevich, Palsson and Sovine~\cite[Propositions~4.1--4.2]{IPS2022},
which builds on the method of Grafakos and Kalton~\cite{GrafakosKalton2001}.

The maximal problem has a different scale dependence. Borges,
Foster, Ou, Pipher and Zhou~\cite{BorgesEtAl2023} obtained sparse
bounds for bilinear spherical maximal averages. The effect of partial
curvature was investigated by Chen, Grafakos, He, Honz\'{\i}k and
Slav\'{\i}kov\'{a}~\cite{ChenEtAl2022}, while Borges, Foster and
Ou~\cite{BorgesFosterOu2025} studied fractal dilation sets. Here we
keep the scale fixed and work directly with the difference projection.

In Section~7.1 of her survey, Tainara Borges~\cite[Question~7.2]{Borges2026}
asks what remains true when the sphere is replaced by a curved hypersurface.
The question asks whether full curvature is enough for the endpoint
\eqref{eq:knownsphere}. The difference map provides a uniform route for general hypersurfaces.
Its critical points determine the projected density, which in turn controls
the strong endpoint.

There is a geometric distinction here that matters for the statement of the
question. Borges works with compact, smooth hypersurfaces without boundary
in~\cite[Section~7.1]{Borges2026}. On each connected component of such a
hypersurface, nonvanishing Gaussian curvature forces the second fundamental
form to be definite. We recall a direct proof below. For further geometric
context, see Ghomi~\cite{Ghomi2002}.

The essential local condition is weaker than global definiteness:
it suffices to have a definite 2D restriction of the second fundamental
form within the kernel of the difference map. In particular, the
hypothesis makes sense on a surface patch. The curved-patch example below exhibits the role of the fiber geometry
in the strong endpoint.

\subsection{The two questions and the main result}\label{subsec:results}
We state the question in the notation used below. A closed hypersurface means a compact embedded hypersurface without boundary. We use Gaussian curvature to mean the product of its principal curvatures. For a smooth hypersurface $S\subseteq\R^d_y\times\R^d_z$ with smooth surface density $\sigma$, let
\begin{equation}\label{eq:operator}
 \A^S_\sigma(f,g)(x)\coloneqq
 \int_S f(x-y)g(x-z)\dd\sigma(y,z),\qquad x\in\R^d.
\end{equation}
The problem concerns this fixed-scale average. For $0<r<1$, the target space is equipped with the strong quasi-norm
\begin{equation}\label{eq:quasinorm}
 \norm h_{L^r(\R^d)}\coloneqq
 \left(\int_{\R^d}|h(x)|^r\dd x\right)^{1/r}.
\end{equation}
In particular, the desired $L^{1/2}$ estimate controls the integral of the square root of the absolute value.

The two questions occur together in Borges's survey~\cite[Section~7.1, pp.~18--19]{Borges2026}. We state the broader curvature problem first. Both statements concern fixed-scale bilinear averages, and all the hypersurfaces below are compact, smooth, and without boundary.

\begin{problem}[Borges, Question 7.1~\cite{Borges2026}]\label{prob:Borges71}
Let $d\ge1$ and $1\le k\le2d-1$. Let $S_k\subseteq\R^{2d}$ be a smooth closed hypersurface with at least $k$ nonzero principal curvatures at each point. Determine the best H\"older bounds
\begin{equation}\label{eq:original71}
 \A^{S_k}_\sigma:L^p(\R^d)\times L^q(\R^d)\longrightarrow L^r(\R^d),
 \qquad \frac1r=\frac1p+\frac1q,\quad 1\le p,q\le\infty.
\end{equation}
\end{problem}

The known sufficient bounds for part of the partially curved range are described in Borges~\cite[Equation~(27)]{Borges2026}, drawing on the hypersurface estimates of Cho, Lee and Shuin~\cite[Proposition~1.2]{ChoLeeShuin2024} and related work of Borges, Foster and Ou~\cite{BorgesFosterOu2025}. Our corollary below gives the full H\"older family on a geometric subclass with partially vanishing principal curvatures.

The second question is the principal objective of this paper.

\begin{problem}[Borges, Question 7.2~\cite{Borges2026}]\label{prob:Borges}
Let $d\ge2$. Suppose $S\subseteq\R^{2d}$ is a compact smooth hypersurface without boundary and its Gaussian curvature never vanishes. Does there exist $C_{S,\sigma}<\infty$ such that
\begin{equation}\label{eq:question}
 \norm{\A^S_\sigma(f,g)}_{L^{1/2}(\R^d)}
 \le C_{S,\sigma}\norm{f}_{L^1(\R^d)}\norm{g}_{L^1(\R^d)}
\end{equation}
for every $f,g\in L^1(\R^d)$ and every smooth surface density $\sigma$ on $S$?
\end{problem}

The constant in this question may depend on the particular surface and on
its density. The known spherical result of Iosevich, Palsson and
Sovine~\cite[Section~5.2]{IPS2022} is therefore a special case. Our first
result treats the whole class in \Probref{prob:Borges}. The second part
shows what can happen when one allows a curved patch outside that class.

\begin{theorem}[Closed hypersurfaces and boundary sensitivity]\label{thm:main}
Let $d\ge2$.
\begin{enumerate}[label=\textup{(\alph*)},leftmargin=*,itemsep=3pt]
\item Let $S\subseteq\R^{2d}$ be a compact smooth embedded hypersurface without boundary. If its Gaussian curvature is nonzero everywhere, then each smooth, possibly signed, density $\sigma$ on $S$ satisfies
\begin{equation}\label{eq:main}
 \norm{\A^S_\sigma(f,g)}_{L^{1/2}(\R^d)}
 \le C_{S,\sigma}\norm{f}_{L^1(\R^d)}\norm{g}_{L^1(\R^d)}.
\end{equation}
\item There are a compact smooth embedded hypersurface patch $S_{\mathrm{patch}}\subseteq\R^{2d}$ with smooth boundary and nowhere-vanishing Gaussian curvature, and a nonnegative smooth density $\sigma$ supported in its interior, for which \eqref{eq:main} fails. More precisely, there are nonnegative $f_n=g_n\in C_c^\infty(\R^d)$ such that
\begin{align}\label{eq:normalizedcounterexample}
 \norm{f_n}_1=\norm{g_n}_1&=1,\\
 \norm{\A^{S_{\mathrm{patch}}}_\sigma(f_n,g_n)}_{1/2}
 &\longrightarrow\infty.\notag
\end{align}
\end{enumerate}
\end{theorem}

Part~\textup{(a)} answers \Probref{prob:Borges} in full. Part~\textup{(b)}
gives a complementary curved-patch obstruction and identifies the geometry
responsible for the contrast. The second fundamental form is definite along
the relevant difference fibers in part~\textup{(a)}, whereas the patch has a
rank-one fiber Hessian at its critical points.

The endpoint in \Thmref{thm:main}\textup{(a)} also gives the complete family of fixed-scale H\"older estimates.
\begin{corollary}[Full H\"older family]\label{cor:holder}
Under the hypotheses of \Thmref{thm:main}\textup{(a)}, let $1\le p,q\le\infty$, where $r^{-1}=p^{-1}+q^{-1}$. Then
\begin{equation}\label{eq:holder}
 \norm{\A^S_\sigma(f,g)}_{L^r(\R^d)}
 \le C_{S,\sigma,p,q}\norm f_{L^p(\R^d)}\norm g_{L^q(\R^d)}.
\end{equation}
The constant remains the same at each fixed dilation $t>0$ when $(y,z)$ is replaced by $(ty,tz)$.
\end{corollary}

The geometry behind the endpoint also gives a useful instance of
\Probref{prob:Borges71}. At a critical point, two definite fiber
directions already suffice for the density estimate.

\begin{corollary}[A partial answer to Problem 1]\label{cor:partial71}
Let $d\ge2$ and let $S\subseteq\R^{2d}$ be a smooth closed hypersurface. Suppose that at every critical point $p$ of $\pn|_S$, the second fundamental form has a definite restriction to some 2D subspace of $\ker\pn$. Then every smooth density satisfies the entire H\"older family \eqref{eq:holder}.

In particular, this conclusion holds when $S$ is the smooth boundary of a compact convex body and at least $d+1$ of its principal curvatures are nonzero at every point.
\end{corollary}

\begin{remark}\label{rem:contribution}
The quasi-Banach localization principle is due to Iosevich, Palsson and
Sovine~\cite[Propositions~4.1--4.2]{IPS2022}. The geometric argument
below verifies its density hypothesis for every closed hypersurface in
\Probref{prob:Borges}. The curved-patch example identifies a rank-one fiber fold as a concrete obstruction to the strong endpoint. The
logarithmic test is related to the 1D circular example of
Wei~\cite[Section~4]{Wei2026}.
\end{remark}

\begin{remark}[A partially flat closed example]\label{rem:flat}
The local criterion also applies to surfaces with flat principal
directions. Write $u=(y+z)/2$ and $v=y-z$ and consider
\begin{equation}\label{eq:flatshort}
 S_{\mathrm{flat}}=\left\{(y,z): |u|^2+\sum_{j=1}^d v_j^4=1\right\}.
\end{equation}
This is a smooth compact hypersurface. At each critical point of
$\pn|_{S_{\mathrm{flat}}}$, we have $u=0$ and the fiber Hessian is
$2I_d$. Thus the strong endpoint follows from
\Lemref{lem:density} and \Lemref{lem:localization}, even though
some principal curvatures vanish. Hence this closed, partially flat
hypersurface has the strong endpoint and the entire H\"older family,
providing a concrete positive example for \Probref{prob:Borges71}.
\end{remark}

\subsection{Organization}\label{subsec:organization}
The main theorem is proved in \Secref{sec:mainproofs}. We state the two
lemmas used in its proof before turning to the positive result and the patch
example. Their proofs appear in \Secref{sec:auxproofs}.

\subsection{Notation}\label{subsec:notation}
We denote the difference map by $\pn(y,z)=y-z$, and the support and total variation of a signed measure by $\spt\mu$ and $|\mu|$. For a Borel map $P$, the measure $P_*\mu$ is defined by $P_*\mu(E)=\mu(P^{-1}(E))$. For nonnegative quantities,
\begin{equation}\label{eq:comparison}
 A\lesssim B\quad\Longleftrightarrow\quad A\le CB,
 \qquad
 A\eqsim B\quad\Longleftrightarrow\quad
 A\lesssim B\ \text{and}\ B\lesssim A.
\end{equation}
Constants may depend on the fixed surface, measure and dimension. When we omit a domain from a Lebesgue norm, that domain is $\R^d$. We begin with nonnegative bounded functions of compact support and extend estimates to general functions by monotone approximation and absolute values.

\section{Proofs of the main results}\label{sec:mainproofs}
\subsection{The analytic and geometric lemmas}\label{subsec:auxstatements}
Write $\pn(y,z)=y-z$. For a compactly supported positive measure $\mu$ on $\R^{2d}$, define
\begin{equation}\label{eq:projintro}
 T_\mu(f,g)(x)=\int_{\R^{2d}}f(x-y)g(x-z)\dd\mu(y,z).
\end{equation}
The first lemma is the analytic passage from a bounded difference
density to the strong endpoint. The second gives a geometric condition
for that density. We state both before proving the main theorem.

\begin{lemma}[Difference-density criterion]\label{lem:localization}
Let $\mu$ be a finite positive Borel measure with $\spt\mu\subseteq[-R,R]^{2d}$ for some $R>0$. Suppose
\begin{equation}\label{eq:differencehyp}
 (\pn)_*\mu=k(v)\dd v,
 \qquad\text{and}\qquad k\in L^\infty(\R^d).
\end{equation}
Then the defining integral is absolutely convergent almost everywhere for $f,g\in L^1(\R^d)$, is independent of their representatives, and satisfies
\begin{equation}\label{eq:localtarget}
 \norm{T_\mu(f,g)}_{1/2}
 \le 9^d R^d\norm k_\infty\norm f_1\norm g_1.
\end{equation}
\end{lemma}
This is the fixed-scale specialization of the method in~\cite[Propositions~4.1--4.2]{IPS2022}. We include its full proof.

\begin{lemma}[Local density criterion for the difference projection]\label{lem:density}
Let $d\ge2$, and let $S\subseteq\R^{2d}$ be a smooth embedded hypersurface, possibly with boundary. Let $a\in C_c^\infty(S^\circ)$ be nonnegative, where $S^\circ$ is the interior of $S$ as a manifold. At every critical point $p\in\spt a$ of $\pn|_S$, assume that
the second fundamental form is definite on a 2D subspace
$E_p\subseteq K$, where
\begin{equation}\label{eq:K}
 K=\ker\pn=\{(h,h):h\in\R^d\}\subseteq T_pS.
\end{equation} Then
\begin{equation}\label{eq:densityresult}
 (\pn)_*(a\dd A_S)=k_a(v)\dd v,\quad\text{and}\quad \norm{k_a}_{L^\infty(\R^d)}<\infty.
\end{equation}
\end{lemma}
For a compact hypersurface without boundary, the choice $a\equiv1$ is allowed. For a patch, the same lemma applies to smooth densities supported away from its boundary.

\begin{lemma}[Curvature of closed hypersurfaces]\label{lem:closeddef}
Let $M\subseteq\R^{n+1}$ be a compact smooth embedded hypersurface without boundary. If all its principal curvatures are nonzero, its second fundamental form is definite at every point.
\end{lemma}

\begin{lemma}[A 1D fold obstruction]\label{lem:obstruction}
If $\theta\in C_c^\infty(\R)$ is nonnegative and positive near $0$, then
\begin{equation}\label{eq:1Dfold}
 \B_\theta(F,G)(x)=\int_\R\theta(t)
 F\left(x-t-\frac{t^2}{2}\right)
 G\left(x-t+\frac{t^2}{2}\right)\dd t
\end{equation}
is unbounded from $L^1(\R)\times L^1(\R)$ to strong $L^{1/2}(\R)$. More precisely, there are $c_1>0$ and nonnegative $F_N\in C_c^\infty(\R)$ with uniformly bounded $L^1$ norms such that
\begin{equation}\label{eq:foldsequence}
 \int_0^{c_1}\B_\theta(F_N,F_N)(x)^{1/2}\dd x\longrightarrow\infty.
\end{equation}
\end{lemma}
The logarithmic test in \Lemref{lem:obstruction} is of the same type as the circular counterexample in Wei~\cite[Section~4]{Wei2026}. Its proof below is given directly for \eqref{eq:1Dfold}.

\subsection{Positive answer to the original question}\label{subsec:positive}
\begin{proof}[Proof of \Thmref{thm:main}\textup{(a)}]
We verify the geometric hypothesis of \Lemref{lem:density}. The map $\pn:\R^{2d}\to\R^d$ is linear, and its kernel is the $d$-plane $K=\{(h,h):h\in\R^d\}$. Since $\dim T_pS=2d-1$, we have
\begin{equation}\label{eq:criticalcriterion}
\begin{gathered}
\dim(T_pS\cap K)\in\{d-1,d\},\\[3pt]
\operatorname{rank}d(\pn|_S)_p=2d-1-\dim(T_pS\cap K).
\end{gathered}
\end{equation}
Therefore $p$ is critical precisely when $K\subseteq T_pS$. At such a point, \Lemref{lem:closeddef} makes the restriction of the second fundamental form to $K$ definite. Since $S$ is compact and has no boundary, the constant function $a\equiv1$ belongs to $C_c^\infty(S)$. We may apply \Lemref{lem:density} with this choice and obtain
\begin{equation}\label{eq:mainpush}
 (\pn)_*\dd A_S=k_S(v)\dd v,\quad\text{and}\quad \norm{k_S}_{L^\infty(\R^d)}<\infty.
\end{equation}

We next pass to the given surface density. Write $\dd\sigma=a\dd A_S$ and choose the positive measure
\begin{equation}\label{eq:majorant}
 M_\sigma=\norm a_{L^\infty(S)},\qquad
 \dd\mu=M_\sigma\dd A_S,\quad\text{and}\quad |\sigma|\le\mu.
\end{equation}
The surface is compact, so $M_\sigma<\infty$. From \eqref{eq:mainpush} we obtain
\begin{align}\label{eq:mainchain}
 (\pn)_*\mu&=M_\sigma k_S(v)\dd v,\\
 |\A^S_\sigma(f,g)(x)|&\le T_\mu(|f|,|g|)(x),\notag\\
 \norm{\A^S_\sigma(f,g)}_{1/2}
 &\le C_{S,\sigma}\norm f_1\norm g_1.\notag
\end{align}
The last line is exactly \Lemref{lem:localization} applied to the preceding two lines.

The absolute convergence and independence of representatives in \Lemref{lem:localization} apply to the positive majorant. They therefore also hold for the original integral against $\sigma$.
\end{proof}

\subsection{A curved patch with a different endpoint}\label{subsec:negative}
\begin{proof}[Proof of \Thmref{thm:main}\textup{(b)}]
We arrange a rank-one quadratic fold in the difference fibers,
while keeping the full second fundamental form nondegenerate.

Use the invertible linear change of variables
\begin{equation}\label{eq:uv}
 u=\frac{y+z}{2},\qquad v=y-z,
 \qquad\text{and}\qquad y=u+\frac v2,\quad z=u-\frac v2.
\end{equation}
Write $v'=(v_2,\ldots,v_d)$ and define the graph
\begin{equation}\label{eq:patchgraph}
 v_1=\Phi(u,v')\coloneqq u_1^2+\sum_{j=2}^d u_jv_j.
\end{equation}
In the variables $(u_1,u_2,v_2,\ldots,u_d,v_d)$, the Hessian consists of the block $(2)$ and $d-1$ blocks $\left(\begin{smallmatrix}0&1\\1&0\end{smallmatrix}\right)$. Consequently,
\begin{equation}\label{eq:hessiandet}
 \det D^2\Phi=2(-1)^{d-1}\neq0.
\end{equation}
The second fundamental form of this graph is a nonzero scalar multiple of $D^2\Phi$. An invertible linear change of ambient coordinates changes the second fundamental form by congruence and a nonzero scalar. Thus \eqref{eq:hessiandet} also gives nonvanishing Gaussian curvature in the original $(y,z)$ coordinates. On the other hand,
\begin{equation}\label{eq:patchfiberrank}
 D_u\Phi=(2u_1,v'),
 \qquad\text{and}\qquad
 D_u^2\Phi=\operatorname{diag}(2,0,\ldots,0).
\end{equation}
At the critical points $u_1=0$ and $v'=0$, the fiber Hessian has rank one. Restrict the graph to a small closed parameter ball and choose a nonnegative density
\begin{equation}\label{eq:patchmeasure}
 \dd\sigma=\rho(u,v')\dd u\dd v',\quad\text{and}\quad \rho\ge c>0
\end{equation}
on a smaller parameter box. The parameter ball gives a smooth boundary, and the density is chosen to vanish near that boundary. Choose a cube $I\subseteq\R^{d-1}$, a nonnegative $\theta\in C_c^\infty(\R)$ positive near zero, and $c>0$ so that
\begin{equation}\label{eq:rhominorant}
 \rho(t,x',0)\ge c\theta(t),\qquad t\in\R,\quad\text{and}\quad x'\in I.
\end{equation}
These choices depend only on the fixed density.

Let $F_N\in C_c^\infty(\R)$ and $c_1>0$ be as in \Lemref{lem:obstruction} for this $\theta$. Take a fixed $h\in C_c^\infty(\R)$ with $h\ge0$ and $\int h=1$. For $\varepsilon>0$ define
\begin{align}\label{eq:tensordata}
 h_\varepsilon(t)&=\varepsilon^{-1}h(t/\varepsilon),\\
 f_{N,\varepsilon}(x)&=g_{N,\varepsilon}(x)
 =F_N(x_1)\prod_{j=2}^d h_\varepsilon(x_j).\notag
\end{align}
Both norms are $\norm{F_N}_1$, independently of $\varepsilon$. We will first keep $N$ fixed and then let $\varepsilon$ tend to zero.

For each $j\ge2$, let
\begin{equation}\label{eq:ab}
 a_j=x_j-u_j-v_j/2,\quad\text{and}\quad b_j=x_j-u_j+v_j/2.
\end{equation}
The Jacobian of $(u_j,v_j)\mapsto(a_j,b_j)$ has absolute value $1$. Hence, after $a_j=\varepsilon s_j$ and $b_j=\varepsilon t_j$,
\begin{align}\label{eq:inverseab}
 u_j&=x_j-\frac{\varepsilon(s_j+t_j)}{2},
 &v_j&=\varepsilon(t_j-s_j),\\
 h_\varepsilon(a_j)h_\varepsilon(b_j)\dd u_j\dd v_j
 &=h(s_j)h(t_j)\dd s_j\dd t_j.\notag
\end{align}
Write $s'=(s_2,\ldots,s_d)$ and $t'=(t_2,\ldots,t_d)$, and define
\begin{align}\label{eq:epsphase}
 u'_\varepsilon&=x'-\frac{\varepsilon(s'+t')}{2},
 &v'_\varepsilon&=\varepsilon(t'-s'),\\
 \Phi_\varepsilon&=u_1^2+
 \varepsilon\sum_{j=2}^d\left(x_j-\frac{\varepsilon(s_j+t_j)}2\right)(t_j-s_j).\notag
\end{align}
The concentration factors cancel the Jacobian. Consequently,
\begin{align}\label{eq:patchexact}
 &\A^S_\sigma(f_{N,\varepsilon},g_{N,\varepsilon})(x)\\
 &=\int\!\int \rho(u_1,u'_\varepsilon,v'_\varepsilon)
 \prod_{j=2}^d h(s_j)h(t_j)
 F_N\left(x_1-u_1-\frac{\Phi_\varepsilon}{2}\right)
 F_N\left(x_1-u_1+\frac{\Phi_\varepsilon}{2}\right)
 \dd s'\dd t'\dd u_1.\notag
\end{align}
For each fixed $N$ and $x$, the integrand is bounded on a compact parameter set. As $\varepsilon\downarrow0$, we have $u'_\varepsilon\to x'$, $v'_\varepsilon\to0$, and $\Phi_\varepsilon\to u_1^2$. Dominated convergence, together with $\int h=1$, gives
\begin{align}\label{eq:patchlimit}
 &\lim_{\varepsilon\downarrow0}
 \A^S_\sigma(f_{N,\varepsilon},g_{N,\varepsilon})(x)\\
 &=\int\rho(u_1,x',0)
 F_N\left(x_1-u_1-\frac{u_1^2}{2}\right)
 F_N\left(x_1-u_1+\frac{u_1^2}{2}\right)\dd u_1.\notag
\end{align}
The choice of density in \eqref{eq:rhominorant} now has its purpose.
For $x'\in I$, the limiting expression is bounded below by
\begin{equation}\label{eq:oneDminorant}
 c\B_\theta(F_N,F_N)(x_1),\qquad x'\in I.
\end{equation}
For fixed $N$, Fatou's lemma and \eqref{eq:patchlimit} give
\begin{equation}\label{eq:patchcontradiction}
\begin{gathered}
\liminf_{\varepsilon\downarrow0}
 \int_{\R^d}\A^S_\sigma(f_{N,\varepsilon},g_{N,\varepsilon})^{1/2}\dd x
 \ge b_N,\\[5pt]
b_N\coloneqq c^{1/2}|I|\int_0^{c_1}
 \B_\theta(F_N,F_N)(x_1)^{1/2}\dd x_1.
\end{gathered}
\end{equation}
By \Lemref{lem:obstruction}, $b_N$ tends to infinity and $m_N\coloneqq\norm{F_N}_1$ is bounded above by a finite constant $M_0$. For all sufficiently large $N$, $m_N>0$. Choose $0<\varepsilon_N<1/N$ so that the integral in \eqref{eq:patchcontradiction} is at least $b_N/2$, and define
\begin{align}\label{eq:diagonalsequence}
 \widetilde f_N=\widetilde g_N
 &\coloneqq m_N^{-1}f_{N,\varepsilon_N},\\
 \norm{\widetilde f_N}_1=\norm{\widetilde g_N}_1&=1,\notag\\
 \norm{\A^S_\sigma(\widetilde f_N,\widetilde g_N)}_{1/2}
 &\ge \left(\frac{b_N}{2M_0}\right)^2
 \longrightarrow\infty.\notag
\end{align}
Thus the smooth functions obtained by this diagonal choice have unit
$L^1$ norms and arbitrarily large $L^{1/2}$ averages. Keeping the two
limits in this order was essential: the passage to $\varepsilon=0$ was
justified only for each fixed $N$. This proves \eqref{eq:normalizedcounterexample}.
\end{proof}

\subsection{The entire H\"older family}\label{subsec:allholder}
\begin{proof}[Proof of \Corref{cor:holder}]
We may take $f,g\ge0$, since absolute values give the general case. Let $\mu=|\sigma|$, $T=T_\mu$ and $m=\mu(S)$. If $m=0$, the assertion is immediate. Assume $m>0$. The smooth area majorant in \eqref{eq:majorant} gives the same endpoint for $T$, even though $|\sigma|$ need not have a smooth density. Together with Tonelli's theorem and $T(1,1)=m$, this gives
\begin{equation}\label{eq:extreme}
\begin{gathered}
\norm{T(F,G)}_{1/2}\le C\norm F_1\norm G_1,\\[3pt]
\norm{T(F,1)}_1=m\norm F_1
\quad\text{and}\quad
\norm{T(1,G)}_1=m\norm G_1.
\end{gathered}
\end{equation}
First assume $1\le p,q<\infty$ and set $\theta=p^{-1}$, $\varphi=q^{-1}$, and $r^{-1}=\theta+\varphi$.
Suppose first that $\theta+\varphi\ge1$, so $r\le1$. Set $a=\theta+\varphi-1$, $b=1-\varphi$, and $c=1-\theta$.
These numbers are nonnegative and satisfy
\begin{equation}\label{eq:exponentcheck}
 a+b+c=1,\quad p(a+b)=q(a+c)=1,\quad\text{and}\quad 2a+b+c=r^{-1}.
\end{equation}
Indeed, at every $(y,z)$ the integrand has the exact factorization
\begin{align}\label{eq:holderfactor}
 f(x-y)g(x-z)
 &=\bigl(f(x-y)^p g(x-z)^q\bigr)^a
 \bigl(f(x-y)^p\bigr)^b
 \bigl(g(x-z)^q\bigr)^c.
\end{align}
H\"older's inequality with exponents $1/a,1/b,1/c$ in the measure variable (omitting a factor if its weight is zero) gives
\begin{equation}\label{eq:interpolate1}
 T(f,g)\le T(f^p,g^q)^a T(f^p,1)^b T(1,g^q)^c.
\end{equation}
For the spatial integral, the relevant exponents are
\begin{equation}\label{eq:spatialexponents}
 \frac1{2ar},\quad\frac1{br},\quad\frac1{cr},
 \qquad\text{and}\qquad 2ar+br+cr=1.
\end{equation}
We omit zero-weight factors as above. These are ordinary H\"older exponents, all at least one. Applying that inequality to the $r$th power of \eqref{eq:interpolate1} and taking the $r$th root, we obtain
\begin{align}\label{eq:interpolate2}
 \norm{T(f,g)}_r
 &\le \norm{T(f^p,g^q)}_{1/2}^{a}
 \norm{T(f^p,1)}_1^b\norm{T(1,g^q)}_1^c\\
 &\le C^a m^{b+c}\norm f_p\norm g_q.\notag
\end{align}
If $\theta+\varphi\le1$, we use the $L^1$ estimates for each marginal. Indeed,
\begin{equation}\label{eq:holderfactor2}
 f(x-y)g(x-z)
 =\bigl(f(x-y)^p\bigr)^\theta
 \bigl(g(x-z)^q\bigr)^\varphi 1^{1-\theta-\varphi}.
\end{equation}
The measure and spatial H\"older inequalities then read
\begin{align}\label{eq:interpolate3}
 T(f,g)&\le T(f^p,1)^\theta
 T(1,g^q)^\varphi m^{1-\theta-\varphi},\\
 \norm{T(f,g)}_r&\le m\norm f_p\norm g_q.\notag
\end{align}
If $p=\infty$, the marginal measure $\nu_z(E)=\mu(\R^d\times E)$ reduces the estimate to convolution, since
\begin{equation}\label{eq:inftycase}
 T(f,g)(x)\le\norm f_\infty(|g|*\nu_z)(x),\quad\text{and}\quad \norm{|g|*\nu_z}_q\le m\norm g_q.
\end{equation}
The case $q=\infty$ is symmetric. This proves \eqref{eq:holder}.

For the fixed-scale dilated average, write $f_t(x)=f(tx)$. We have
\begin{equation}\label{eq:dilation}
\begin{gathered}
\A^S_{\sigma,t}(f,g)(tx)=\A^S_{\sigma,1}(f_t,g_t)(x),\\[3pt]
\norm{f_t}_p=t^{-d/p}\norm f_p
\quad\text{and}\quad
\norm{g_t}_q=t^{-d/q}\norm g_q.
\end{gathered}
\end{equation}
As $d/r=d/p+d/q$, the scale factors cancel. The bound has the same constant for all fixed $t>0$.
\end{proof}

\subsection{The partial-curvature consequence}\label{subsec:partial71}
\begin{proof}[Proof of \Corref{cor:partial71}]
At every critical point $p$ we have $K=\ker\pn\subseteq T_pS$. The first assertion is therefore \Lemref{lem:density}, followed by \Lemref{lem:localization} and the H\"older argument above.

Suppose now that $S$ is convex and has at least $d+1$ nonzero principal curvatures. Choose the normal so that the second fundamental form is nonnegative, and denote by $P_p\subseteq T_pS$ the span of $d+1$ strictly positive principal directions. Since $\dim T_pS=2d-1$ and $\dim K=d$,
\begin{equation}\label{eq:partial71dimension}
 \dim(P_p\cap K)\ge(d+1)+d-(2d-1)=2.
\end{equation}
The second fundamental form is positive definite on any 2D subspace of $P_p\cap K$. This verifies the hypothesis at every critical point, and the stated H\"older bounds follow.
\end{proof}

\section{Proofs of the auxiliary lemmas}\label{sec:auxproofs}
\subsection{From a bounded difference density to the strong endpoint}\label{subsec:locproof}
\begin{proof}[Proof of \Lemref{lem:localization}]
Take $f,g\ge0$ bounded and compactly supported. With $u=x-y$ and $v=y-z$, we have $x-z=u+v$. Tonelli's theorem gives
\begin{align}\label{eq:TlocalL1}
 \int T_\mu(f,g)(x)\dd x
 &=\int_{\R^d}k(v)\int_{\R^d}f(u)g(u+v)\dd u\dd v\\
 &\le\norm k_\infty\norm f_1\norm g_1.\notag
\end{align}
We use cubes at the scale of $\spt\mu$. For $j\in\mathbb Z^d$ set
\begin{align}\label{eq:cubes}
 Q_j&=R\bigl(j+[0,1)^d\bigr),
 &Q_j^*&=R\bigl(j+[-1,2)^d\bigr),\\
 a_j&=\int_{Q_j^*}f,
 &b_j&=\int_{Q_j^*}g.\notag
\end{align}
If $x\in Q_j$ and $(y,z)\in\spt\mu$, then $x-y,x-z\in Q_j^*$. Hence
\begin{equation}\label{eq:finiteinteract}
 T_\mu(f,g)(x)
 =T_\mu(f\one_{Q_j^*},g\one_{Q_j^*})(x),\qquad x\in Q_j.
\end{equation}
Cauchy--Schwarz on $Q_j$, followed by \eqref{eq:TlocalL1}, now yields
\begin{align}\label{eq:cubehalf}
 \int_{Q_j}T_\mu(f,g)^{1/2}\dd x
 &\le R^{d/2}
 \left(\int_{Q_j}T_\mu(f\one_{Q_j^*},g\one_{Q_j^*})\dd x\right)^{1/2}\\
 &\le R^{d/2}\norm k_\infty^{1/2}(a_jb_j)^{1/2}.\notag
\end{align}
The enlarged cubes overlap exactly $3^d$ times almost everywhere, so
\begin{equation}\label{eq:overlap}
 \sum_j\one_{Q_j^*}=3^d,
 \qquad
 \sum_j a_j=3^d\norm f_1,
 \qquad\text{and}\qquad
 \sum_j b_j=3^d\norm g_1.
\end{equation}
We sum \eqref{eq:cubehalf} and apply discrete Cauchy--Schwarz:
\begin{align}\label{eq:summing}
 \int_{\R^d}T_\mu(f,g)^{1/2}\dd x
 &\le R^{d/2}\norm k_\infty^{1/2}\sum_j(a_jb_j)^{1/2}\\
 &\le R^{d/2}\norm k_\infty^{1/2}
   \left(\sum_j a_j\right)^{1/2}
   \left(\sum_j b_j\right)^{1/2}\notag\\
 &=3^dR^{d/2}\norm k_\infty^{1/2}\norm f_1^{1/2}\norm g_1^{1/2}.\notag
\end{align}
Squaring yields \eqref{eq:localtarget}, independently of the number of cubes.

For arbitrary nonnegative $f,g\in L^1$, choose finite Borel representatives and define
\begin{equation}\label{eq:monotoneextension}
f_n=\min\{f,n\}\one_{[-n,n]^d}
\quad\text{and}\quad
g_n=\min\{g,n\}\one_{[-n,n]^d}.
\end{equation}
These sequences are nondecreasing and converge pointwise to $f$ and $g$. Monotone convergence, first in the measure variable and then in $x$, proves the same bound. In particular, $T_\mu(f,g)$ is finite almost everywhere. Absolute values give the result for signed or complex-valued functions.

Finally, let $E\subseteq\R^d$ have Lebesgue measure zero. We have
\begin{equation}\label{eq:representatives}
 \int_{\R^d}\int_{\R^{2d}}\one_E(x-y)\dd\mu(y,z)\dd x
 =\mu(\R^{2d})|E|=0,
\end{equation}
and the same identity holds with $x-z$. Thus changing either function on a null set leaves the average unchanged for almost every $x$. This also justifies working with finite Borel representatives throughout.
\end{proof}

\subsection{The fold density of the difference projection}\label{subsec:densityproof}
\begin{proof}[Proof of \Lemref{lem:density}]
Choose a finite smooth partition of unity on $\spt a$, with chart
supports contained in $S^\circ$. We write $a$ for one localized
density and work in the coordinates
\begin{equation}\label{eq:densitycoords}
 u=\frac{y+z}{2},\qquad v=y-z,\quad\text{and}\quad \pn(u,v)=v.
\end{equation}
Since $\dim T_pS=2d-1$ and $\dim\ker\pn=d$,
\begin{equation}\label{eq:densityrank}
 \ker d(\pn|_S)_p=T_pS\cap\ker\pn,\quad\text{and}\quad \operatorname{rank}d(\pn|_S)_p\in\{d-1,d\}.
\end{equation}
At regular points, the submersion theorem gives coordinates $(t,v)$,
$t\in\R^{d-1}$, in which
\begin{align}\label{eq:regularpush}
 a\dd A_S&=b(t,v)\dd t\dd v,\\
 (\pn)_*(a\dd A_S)&=\left(\int b(t,v)\dd t\right)\dd v.\notag
\end{align}
The chart density has compact support in the $(t,v)$ variables. Its integral in $t$ is therefore a smooth, bounded function of $v$.

Let $p$ be critical. Then $K\subseteq T_pS$. Split the $v$
coordinates as $v=(w,s)$, $w\in\R^{d-1}$, choosing $s$ so that
its direction is transverse to $T_pS$. Locally,
\begin{equation}\label{eq:graphnormal}
 s=\Phi(u,w),\qquad
 \pn|_S(u,w)=(w,\Phi(u,w)),\quad\text{and}\quad D_u\Phi(p)=0.
\end{equation}
For $L(u,v)=(u+v/2,u-v/2)$ write
$X(u,w)=L(u,(w,\Phi(u,w)))$. If $\nu_p$ is a unit normal,
then the second fundamental form in the fiber directions is
\begin{align}\label{eq:IIfiber}
 \partial_{u_i u_j}X(p)&=\partial_{u_i u_j}\Phi(p)\,L(0,e_s),\\
 \mathrm{II}_p\bigl(L(h,0),L(h,0)\bigr)
 &=c_p h^{\mathsf T}D_u^2\Phi(p)h,\notag\\
 c_p&=-\langle\nu_p,L(0,e_s)\rangle\ne0.\notag
\end{align}
The hypothesis therefore supplies two fiber coordinates
$z\in\R^2$ in which the Hessian $D_z^2\Phi(p)$ is definite. Write
$u=(z,t)$, $t\in\R^{d-2}$. If $d=2$, the $t$ variable is absent.

Treat $(t,w)$ as parameters. Solving $D_z\Phi=0$ gives a smooth
critical branch $z_*(t,w)$ and its critical value $s_*(t,w)$:
\begin{equation}\label{eq:branch}
 D_z\Phi(z_*(t,w),t,w)=0,\quad\text{and}\quad s_*(t,w)=\Phi(z_*(t,w),t,w).
\end{equation}
With $h=z-z_*(t,w)$, the exact Taylor formula reads
\begin{align}\label{eq:H}
 \Phi(z,t,w)-s_*(t,w)&=h^{\mathsf T}H(z,t,w)h,\\
 H(z,t,w)&=\int_0^1(1-\tau)
 D_z^2\Phi(z_*(t,w)+\tau h,t,w)\dd\tau.\notag
\end{align}
After shrinking the chart, $\varepsilon H$ is positive definite
for a fixed $\varepsilon\in\{-1,1\}$. Set
\begin{equation}\label{eq:quadratic}
 \eta=(\varepsilon H(z,t,w))^{1/2}h.
\end{equation}
The positive matrix square root is smooth. At $h=0$,
\begin{equation}\label{eq:etaJac}
 D_z\eta=(\varepsilon D_z^2\Phi(z_*(t,w),t,w)/2)^{1/2},\quad\text{and}\quad \det D_z\eta\ne0.
\end{equation}
The inverse function theorem now gives an exact local normal form
\begin{equation}\label{eq:normalform}
 \pn|_S(\eta,t,w)
 =\bigl(w,s_*(t,w)+\varepsilon|\eta|^2\bigr).
\end{equation}

The localized area measure becomes
$b(\eta,t,w)\dd\eta\dd t\dd w$, with smooth compactly supported
$b\ge0$. Indeed, if $J$ is the area Jacobian of $X$, then
\begin{equation}\label{eq:chartdensity}
 b(\eta,t,w)=a(X(z,t,w))J(z,t,w)
 \bigl|\det D_z\eta(z,t,w)\bigr|^{-1}.
\end{equation}
For a nonnegative test function $\varphi$, polar coordinates
$\eta=r\omega$, $\omega\in\mathbb S^1$, give
\begin{align}\label{eq:polartest}
 &\int\varphi(\pn(\eta,t,w))b(\eta,t,w)
 \dd\eta\dd t\dd w\\
 &=\int\!\int\!\int_0^\infty\!\int_{\mathbb S^1}
 b(r\omega,t,w)
 \varphi(w,s_*(t,w)+\varepsilon r^2)
 r\dd\omega\dd r\dd t\dd w.\notag
\end{align}
For fixed $(t,w)$, set $s=s_*(t,w)+\varepsilon r^2$.
The radial Jacobian satisfies $|\dd s|=2r\dd r$.
For $\varepsilon(s-s_*(t,w))>0$, we write
$r=\sqrt{\varepsilon(s-s_*(t,w))}$ and obtain the density
\begin{equation}\label{eq:foldformula}
 k_{\mathrm{loc}}(w,s)
 =\frac12\int_{\R^{d-2}}
 \mathbf1_{\{\varepsilon(s-s_*(t,w))>0\}}
 \int_{\mathbb S^1} b(r\omega,t,w)\dd\omega\dd t.
\end{equation}
The integrand is set to zero outside the local coordinate domain.
Let $B$ contain the projection of $\spt b$ onto the $t$ variables.
When $d=2$, the $t$ integral is omitted and $|B|=1$. Thus
\begin{equation}\label{eq:folddensitybound}
 0\le k_{\mathrm{loc}}(w,s)
 \le \frac{|\mathbb S^1|}{2}|B|\norm b_\infty<\infty.
\end{equation}
In two fiber directions, the polar factor $r$ cancels the
factor $2r$ in $|\dd s|$. In particular, the bound is uniform
as $s$ approaches the critical value $s_*(t,w)$.
For every fixed $(t,w)$, the set $\{\eta=0\}$ has zero
$\dd\eta$ measure. Integrating this statement against
$b(\eta,t,w)\dd t\dd w$ shows that the critical locus has zero
source measure. The polar change of variables therefore accounts
for the entire pushforward, without an extra singular measure. Together with \eqref{eq:regularpush}, the
finite partition gives
\begin{equation}\label{eq:densityglobal}
 (\pn)_*(a\dd A_S)=k_a(v)\dd v,\quad\text{and}\quad \norm{k_a}_\infty<\infty.
\end{equation}
\end{proof}

\subsection{Definiteness on a closed hypersurface}\label{subsec:globaldefproof}
\begin{proof}[Proof of \Lemref{lem:closeddef}]
Fix one connected component $M_0$. Choose $a\notin M_0$ and let $p\in M_0$ maximize $|x-a|^2$. Write
\begin{equation}\label{eq:maxdist}
 R=|p-a|>0,\quad\text{and}\quad \nu=\frac{p-a}{R}.
\end{equation}
The first derivative test makes $\nu$ a unit normal at $p$. For $X\in T_pM_0$, take a smooth curve $\gamma$ on $M_0$ with $\gamma(0)=p$ and $\gamma'(0)=X$. At the maximum,
\begin{align}\label{eq:curvaturemax}
 0&\ge\left.\frac{\dd^2}{\dd t^2}|\gamma(t)-a|^2\right|_{t=0}\\
 &=2|X|^2+2R\langle\nu,\gamma''(0)\rangle.\notag
\end{align}
With $\mathrm{II}_p(X,X)=-\langle\nu,\gamma''(0)\rangle$, it follows that
\begin{equation}\label{eq:positiveatp}
 \mathrm{II}_p(X,X)\ge R^{-1}|X|^2.
\end{equation}
Thus the second fundamental form is definite at $p$. At another
point of the same component choose a continuous local unit normal.
The inertia of its symmetric second fundamental form is locally
constant because no principal curvature vanishes. Reversing the
local unit normal reverses all eigenvalue signs, while preserving
definiteness. Thus the set of points where the second fundamental
form is definite is both open and closed in $M_0$. It contains $p$,
and connectedness shows that it is all of $M_0$. This argument
applies to every component.
\end{proof}

\subsection{The logarithmic test in one dimension}\label{subsec:foldtest}
\begin{proof}[Proof of \Lemref{lem:obstruction}]
Choose $0<c_0<1$ so that $\theta(t)\ge c>0$ for $|t|\le4c_0$. Set
\begin{equation}\label{eq:singularF}
 F(s)=\frac{\one_{(0,c_0)}(s)}{s[\log(e/s)]^{5/4}}.
\end{equation}
This is an $L^1$ function. Indeed, for $0<\delta<c_0$,
\begin{equation}\label{eq:integralF}
 \int_0^\delta F(s)\dd s
 =\frac{4}{[\log(e/\delta)]^{1/4}},\quad\text{and}\quad \norm F_1<\infty.
\end{equation}
The singularity of $F$ is integrable, but just barely. The second
factor in the average will be sampled at distance comparable to $x^2$
from the origin. This produces the extra logarithmic divergence. To make
it precise, fix $0<x<c_1$ and use
\begin{equation}\label{eq:solvet}
 s=x-t-\frac{t^2}{2},\quad\text{and}\quad t=-1+\sqrt{1+2(x-s)}.
\end{equation}
On $0<s<x^2/100$, the relevant branch has $t>0$. After reducing $c_1$ if necessary, direct substitution in \eqref{eq:solvet} gives
\begin{align}\label{eq:foldbounds}
 \frac{x}{2}&\le t\le2x,
 &\frac14x^2&\le s+t^2\le5x^2,\\
 \frac12&\le\left|\frac{\dd t}{\dd s}\right|=\frac1{1+t}\le2,
 &\theta(t)&\ge c.\notag
\end{align}
The two arguments of $F$ in \eqref{eq:1Dfold} are $s$ and $s+t^2$. For $x$ small the latter lies in $(0,c_0)$, and the bounds in \eqref{eq:foldbounds} imply
\begin{equation}\label{eq:secondfactor}
 F(s+t^2)\gtrsim
 \frac1{x^2[\log(e/x)]^{5/4}}.
\end{equation}
Consequently,
\begin{align}\label{eq:lowerB}
 \B_\theta(F,F)(x)
 &\gtrsim\int_0^{x^2/100}F(s)F(s+t(s)^2)\dd s\\
 &\gtrsim\frac{1}{x^2[\log(e/x)]^{5/4}}
   \int_0^{x^2/100}F(s)\dd s\notag\\
 &=\frac{4}{x^2[\log(e/x)]^{5/4}
   [\log(100e/x^2)]^{1/4}}\notag\\
 &\gtrsim\frac1{x^2[\log(e/x)]^{3/2}}.\notag
\end{align}
The power of the logarithm is now below the integrability threshold.
After taking square roots, we obtain
\begin{equation}\label{eq:divergence}
 \int_0^{c_1}\B_\theta(F,F)(x)^{1/2}\dd x
 \gtrsim\int_0^{c_1}\frac{\dd x}{x[\log(e/x)]^{3/4}}
 =\infty.
\end{equation}

Finally, we replace $F$ by smooth functions without losing the contradiction. Take $\chi\in C_c^\infty((-c_0,c_0))$, with $0\le\chi\le1$ and $\chi=1$ on $(0,c_0/2)$. Let $\zeta$ be a smooth nondecreasing function with $\zeta(t)=0$ for $t\le1/2$ and $\zeta(t)=1$ for $t\ge1$. Define
\begin{equation}\label{eq:FN}
 F_N(s)=\chi(s)F(s)\zeta(Ns)\in C_c^\infty(\R).
\end{equation}
The sequence $F_N$ is nondecreasing and converges pointwise to $\chi F$. Therefore
\begin{align}\label{eq:FNdiv}
 \sup_N\norm{F_N}_1&\le\norm F_1,\\
 \lim_{N\to\infty}\int_0^{c_1}
 \B_\theta(F_N,F_N)(x)^{1/2}\dd x&=\infty.\notag
\end{align}
The second line follows from \eqref{eq:divergence} and monotone convergence, since the region used in the lower bound lies where $\chi=1$. Thus no uniform strong $L^{1/2}$ estimate can hold.
\end{proof}

\medskip
\noindent{\bf Funding }
Not applicable.
\medskip

\noindent{\bf Data Availability} This manuscript has no associated data.

\medskip
\noindent{\bf\Large Declarations}
\medskip

\noindent\textbf{Conflicts of Interest}\space The author declares no conflicts of interest.

\medskip
\noindent{\bf AI disclosure}
During the preparation of this manuscript, the author used ChatGPT and Codex to discuss possible proof strategies, assist in drafting and checking parts of the mathematical arguments, organize references, and improve the clarity and readability of the exposition. All mathematical statements and proofs in the final manuscript were written, reviewed, and verified by the author, who takes full responsibility for the content of this manuscript.


\begin{thebibliography}{99}
\bibitem{BorgesEtAl2023}
T.~Borges, B.~Foster, Y.~Ou, J.~Pipher, Z.~Zhou,
Sparse bounds for the bilinear spherical maximal function,
J. Lond. Math. Soc. (2) 107 (4) (2023) 1409--1449.
\href{https://doi.org/10.1112/jlms.12715}{doi:10.1112/jlms.12715}.

\bibitem{BorgesFosterOu2025}
T.~Borges, B.~Foster, Y.~Ou,
Sobolev smoothing estimates for bilinear maximal operators with fractal dilation sets,
J. Funct. Anal. 288 (2025) 110694.
\href{https://doi.org/10.1016/j.jfa.2024.110694}{doi:10.1016/j.jfa.2024.110694}.

\bibitem{ChoLeeShuin2024}
C.-h.~Cho, J.~B.~Lee, K.~Shuin,
$L^p$ improving properties and maximal estimates for certain multilinear averaging operators,
Rev. Mat. Iberoam. 40 (5) (2024) 1799--1832.
\href{https://doi.org/10.4171/RMI/1479}{doi:10.4171/RMI/1479}.

\bibitem{Borges2026}
T.~Borges,
Survey on bilinear spherical averages and associated maximal operators,
J. Geom. Anal. 36 (2026) Article 94.
\href{https://doi.org/10.1007/s12220-026-02336-w}{doi:10.1007/s12220-026-02336-w}.

\bibitem{ChenEtAl2022}
J.~Chen, L.~Grafakos, D.~He, P.~Honz\'{\i}k, L.~Slav\'{\i}kov\'{a},
Bilinear maximal functions associated with surfaces,
Proc. Amer. Math. Soc. 150 (4) (2022) 1635--1639.
\href{https://doi.org/10.1090/proc/15792}{doi:10.1090/proc/15792}.

\bibitem{Ghomi2002}
M.~Ghomi,
Gauss map, topology, and convexity of hypersurfaces with nonvanishing curvature,
Topology 41 (1) (2002) 107--117.
\href{https://doi.org/10.1016/S0040-9383(00)00028-8}{doi:10.1016/S0040-9383(00)00028-8}.

\bibitem{GrafakosKalton2001}
L.~Grafakos, N.~Kalton,
Some remarks on multilinear maps and interpolation,
Math. Ann. 319 (1) (2001) 151--180.
\href{https://doi.org/10.1007/PL00004426}{doi:10.1007/PL00004426}.

\bibitem{IPS2022}
A.~Iosevich, E.~A.~Palsson, S.~R.~Sovine,
Simplex averaging operators: quasi-Banach and $L^p$-improving bounds in lower dimensions,
J. Geom. Anal. 32 (3) (2022) Article 87.
\href{https://doi.org/10.1007/s12220-021-00843-6}{doi:10.1007/s12220-021-00843-6}.

\bibitem{JeongLee2020}
E.~Jeong, S.~Lee,
Maximal estimates for the bilinear spherical averages and the bilinear Bochner--Riesz operators,
J. Funct. Anal. 279 (7) (2020) 108629.
\href{https://doi.org/10.1016/j.jfa.2020.108629}{doi:10.1016/j.jfa.2020.108629}.

\bibitem{Oberlin1988}
D.~M.~Oberlin,
Multilinear convolutions defined by measures on spheres,
Trans. Amer. Math. Soc. 310 (2) (1988) 821--835.
\href{https://doi.org/10.2307/2000993}{doi:10.2307/2000993}.

\bibitem{Wei2026}
X.~Wei,
End point estimates for the bilinear circular averaging operator,
\arxiv{2610.24465} (2026).

\end{thebibliography}
\end{document}